\documentclass[12pt]{amsart}
\usepackage{amssymb, enumerate}

\theoremstyle{plain}

\newtheorem{thm}{Theorem}[section]
\newtheorem{cor}[thm]{Corollary}
\newtheorem{lem}[thm]{Lemma}
\newtheorem{prop}[thm]{Proposition}

\theoremstyle{definition}
\newtheorem{defi}[thm]{Definition}
\newtheorem{conj}[thm]{Conjecture}
\newtheorem{conv}[thm]{Convention}
\newtheorem{nota}[thm]{Notation}
\newtheorem{rem}[thm]{Remark}
\newtheorem{rems}[thm]{Remarks}
\newtheorem{exa}[thm]{Example}
\newtheorem{exas}[thm]{Examples}
\newtheorem{sit}[thm]{}

\newcommand{\brem}{\begin{rem}}
\newcommand{\brems}{\begin{rems}}
\newcommand{\erem}{\end{rem}}
\newcommand{\erems}{\end{rems}}
\newcommand{\bexa}{\begin{exa}}
\newcommand{\bexas}{\begin{exas}}
\newcommand{\eexa}{\end{exa}}
\newcommand{\eexas}{\end{exas}}
\newcommand{\bdefi}{\begin{defi}}
\newcommand{\edefi}{\end{defi}}
\newcommand{\bcor}{\begin{cor}}
\newcommand{\ecor}{\end{cor}}
\newcommand{\blem}{\begin{lem}}
\newcommand{\elem}{\end{lem}}
\newcommand{\bconv}{\begin{conv}}
\newcommand{\econv}{\end{conv}}
\newcommand{\bconj}{\begin{conj}}
\newcommand{\econj}{\end{conj}}
\newcommand{\bprop}{\begin{prop}}
\newcommand{\eprop}{\end{prop}}
\newcommand{\bthm}{\begin{thm}}
\newcommand{\ethm}{\end{thm}}
\newcommand{\bnota}{\begin{nota}}
\newcommand{\enota}{\end{nota}}
\newcommand{\bsit}{\begin{sit}}
\newcommand{\esit}{\end{sit}}
\newcommand{\be}{\begin{eqnarray}}
\newcommand{\ee}{\end{eqnarray}}

\def\ba{\begin{array}}
\def\ea{\end{array}}
\def\bma{\begin{matrix}}
\def\ema{\end{matrix}}
\def\ben{\begin{enumerate}}
\def\een{\end{enumerate}}

\newcommand{\Spec}{\operatorname{Spec}}
\newcommand{\Frac}{\operatorname{Frac}}

\newcommand{\Der}{\operatorname{Der}}

\newcommand{\Pic}{{\operatorname{Pic}}}

\renewcommand{\div}{{\operatorname{div}}}

\newcommand{\lcm}{{\operatorname{lcm}}}
\newcommand{\rk}{{\operatorname{rk}}}
\newcommand{\supp}{{\operatorname{supp}}}

\newcommand{\A}{{\mathbb A}}

\newcommand{\C}{{\mathbb C}}
\newcommand{\Q}{{\mathbb Q}}
\newcommand{\Z}{{\mathbb Z}}

\newcommand{\p}{{\partial}}

\newcommand{\tf}{{\tilde f}}
\newcommand{\tX}{{\tilde X}}

\title{On a result of Miyanishi-Masuda}

\author{Hubert Flenner}
\address{Fakult\"at f\"ur Mathematik,
Ruhr Universit\"at Bochum,
Geb.\ NA 2/72,
Universit\"ats\-stra\ss e\ 150,
44780 Bochum, Germany}
\email{Hubert.Flenner@rub.de}

\author{Mikhail Zaidenberg}
\address{Universit\'e
Grenoble I, Institut Fourier, UMR 5582 CNRS-UJF, BP 74,
38402 St.\ Martin
d'H\`eres c\'edex, France}
\email{zaidenbe@ujf-grenoble.fr}

\thanks{ This research was done during a visit of the first
author at the Institut Fourier of
the University of Grenoble.
He thanks this institution for
its hospitality and  support.\\
\mbox{\hspace{11pt}}{\it 2000 Mathematics Subject Classification}:
13A02, 13F15, 14R05, 14L30.\\
\mbox{\hspace{11pt}}{\it Key words}: affine pseudo-plane,
$\C^*$-action, $\C_+$-action, affine surface}

\begin{document}

\maketitle

\section{Introduction}

Let $X$ be a smooth affine surface over $\C$ with an affine ruling
(an $\A^1$-fibration) $\rho:X\to \A^1_\C$. Assume that $\rho$ is
surjective, has a unique degenerate fiber, and this fiber is
irreducible. In \cite{MaMi} such a surface $X$ is called {\it
affine pseudo-plane}. It is {\it of class} ML$_1$ if $\rho$ is
essentially unique, that is for any other affine ruling
$\rho':X\to \A^1_\C$, the general fibers of $\rho$ and $\rho'$
are the same. In \cite{MaMi} the following classification result
is obtained.

\bthm\label{MT} {\em (Miyanishi-Masuda)} Suppose that $X$ is an
affine pseudo-plane of class ML$_1$. If $X$ admits an effective
$\C^*$-action then the following hold.

\begin{enumerate} [(i)] \item  This $\C^*$-action is
necessarily hyperbolic. \item The universal covering $\tf:\tX\to
X$ is a cyclic covering of degree $d$, where $d$ is the
multiplicity of the unique degenerate fiber of $\rho$. \item
$\tX$ is an affine hypersurface in $\A^3_\C=\Spec\C[x,y,z]$ with
equation $x^m y=z^d-1$ for some $m> 1$. \item The Galois group
$\Z_d=\langle\zeta\rangle$ of the covering $\tf:\tX\to X$, where
$\zeta=\zeta_d$ is a primitive $d$-th root of unity, acts on
$\tX$ via $\zeta.(x,y,z)=(\zeta x, \zeta^{-m}y, \zeta^e z)$, where
$\gcd (e,d)=1$.
\item The $\C^*$-action $\lambda.(x,y,z):=(\lambda x,
\lambda^{-m}y,z)$ ($\lambda\in \C^*$) on $\tX$ descends to the
given $\C^*$-action on $X$,  up to replacing $\lambda$ by
$\lambda^{-1}$.
\end{enumerate}\ethm

Let us add some remarks. An affine ruling on $X$ induces an
affine ruling $\tilde\rho:\tX\to\A^1_\C$ with a unique degenerate
fiber consisting of $d$ disjoint components isomorphic to
$\A^1_\C$. In case $m>1$ there is an essentially unique such
affine ruling on $\tX$, defined by the restriction $x\vert \tX$.
However, for $m=1$, $y\vert \tX$ gives a second independent
affine ruling, which also descends to $X=\tX/\Z_d$. Thus in this
case $X$ cannot be a ML$_1$ surface.

If we want the $\Z_d$-action on $\tX$ to be free, the exponents
$e$ and $d$ above must be coprime. Indeed, otherwise
$\zeta^{eb}=1$ for some $b$ with $0<b<d$, and we would have
$\zeta^b.(0,0,z)=(0,0,z)$ for every $d$-th root of unity $z$.

On the other hand, for every triple $(d,e,m)$ with $d\ge 1, m\ge
2$ and  $\gcd (e,d)=1$, (iii)-(v) determine a smooth affine
pseudo-plane $X$ of class ML$_1$ with an effective $\C^*$-action.
Thus Theorem \ref{MT} provides indeed a complete classification
of these surfaces.

Here we give an alternative proof of Theorem \ref{MT} based on the
results in \cite{FlZa1, FlZa2}.

\section{The proof}

\noindent Under the assumptions of Theorem \ref{MT} $X\not\cong
\A^2_\C$, since otherwise $X$ would admit another affine ruling
$\rho':X\to \A^1_\C$ with general fibers different from those of
$\rho$, which contradicts the condition ML$_1$.

A smooth affine surface $X$ with an elliptic $\C^*$-action is
always isomorphic to $\A^2_\C$, so this case is impossible. If
$X$ is smooth and the $\C^*$-action on $X$ is parabolic then
according to Proposition 3.8(b) in \cite{FlZa1}, $X=\Spec A_0[D]$
for an integral divisor $D$ on a smooth affine curve $C=\Spec
A_0$. The existence of an affine ruling $\rho$  on $X$ with the
base $\A^1_\C$ implies that $C\cong \A^1_\C$. Hence $D$ is a
principal divisor. By Theorem 3.2(b) in \cite{FlZa1},  we have
again $X\cong \A^2_\C= \Spec A_0[0]$ with $A_0=\C[t]$, which
is impossible.

Thus the $\C^*$-action on $X=\Spec A$ is necessarily hyperbolic.
Accordingly we can write \be\label{E1} A=A_0[D_+,D_-] \ee with a
pair of $\Q$-divisors $D_\pm$ on a smooth affine curve $C=\Spec
A_0$ satisfying $D_++D_-\le 0$, see Theorem 4.3 in \cite{FlZa1}.
The remainder of the proof is based on Lemmas \ref{L1} and
\ref{L2} below.

\blem\label{L1} Under the assumptions of Theorem \ref{MT},
$A\cong A_0[D_+,D_-]$, where $A_0=\C[t]$ and $$D_+=-\frac{e'}{d}
[0],\qquad D_-=\frac{e'}{d} [0] -\frac{1}{m} [1]\,.$$ \elem

\smallskip

\noindent {\it Proof of Lemma \ref{L1}.} By Lemmas 1.6 and 2.1 in
\cite{FlZa2}, $X$ admits an affine ruling over an affine base if
and only if it admits a non-trivial $\C_+$-action defined by a
non-zero homogeneous locally nilpotent derivation $\partial\in
\Der (A)$. Moreover, $A_0=\C[t]$ in (\ref{E1}) and, up to an
automorphism $\lambda\longmapsto\lambda^{-1}$ of $\C^*$ (thus
switching $(D_+,D_-)\longmapsto (D_-,D_+)$) we may assume that
$e=\deg
\partial \ge 0$. By Lemma 3.5 and Corollary 3.27 in \cite{FlZa2},
$e=0$ implies that $X\cong \A^1_\C\times \C^*$, so the induced
affine ruling $X\to \C^*$ is essentially unique and has the base
$\C^*$, which contradicts our assumption. Thus $e>0$.

According to Corollary 3.23 in \cite{FlZa2}, the latter implies
that the fractional part $\{D_+\}=D-\lfloor D\rfloor$ is zero or
is supported on one point, and we can choose this point to be
$0\in\A^1_\C$. Such a surface $X=\Spec A$ is of class ML$_1$ if
and only if the fractional part $\{D_-\}$ is supported on at
least 2 points, see \cite[Theorem 4.5]{FlZa2}.

Replacing $(D_+,D_-)$ by the equivalent pair
$(\{D_+\},\,D_-+\lfloor D_+\rfloor)$ (see Theorem 4.3(b) in
\cite{FlZa1}) we may suppose that $D_+=\{D_+\}=-e'/d [0]$, where
$\gcd (e',d)=1$ and $d>0$.

For any affine pseudo-plane $X$, the Picard group $\Pic X$ is a
torsion group \cite[Ch. 3, 2.4.4]{Mi}. On the other hand, for a
$\C^*$-surface $X$ as above, $\rk_\Q (\Pic X\otimes\Q)\ge l-1$,
where $l$ is the number of points $b_j\in\A^1_\C$ such that
$(D_++D_-)(b_j)<0$, see Corollary 4.24 in \cite{FlZa2}. Hence
$l\le 1$ and so, $\exists p\in\A^1_\C : (D_++D_-)(q)=0$ $\forall
q\neq p$.

Since $D_+(q)=0$ $\forall q\neq 0$ we have $D_-(q)=0$ $\forall
q\neq 0,p$. It follows that $\supp (D_-)=\supp (\{D_-\})=\{0,p\}$
with $p\neq 0$. After an automorphism of $\A^1_\C$ we may assume
that $p=1$.
Thus finally $$D_\pm (0)=\mp e'/d,\qquad D_+(1)=0,\quad
D_-(1)=a/m\not\in\Z\qquad  \mbox{and}\quad D_\pm
(q)=0\,\,\,\forall q\neq 0,1\,,$$ where $\gcd (a,m)=1$ and $m>0$.
The smoothness of $X$ forces $a=-1$, see Theorem 4.15 in
\cite{FlZa1}. This proves Lemma \ref{L1}. \qed

Next we use the following description \cite[Corollary
3.30]{FlZa2}, where for a $\Q$-divisor $D$, $d(D)$ denotes the
minimal positive integer $d$ such that $dD$ is integral.

\blem\label{L2} We let $A=\C[t][D_+,D_-]$, where $D_++D_-\le 0$,
$d(D_+)=d$, $d(D_-)=k$. We assume that $D_+=-\frac{e'}{d}[0]$ and
$D_-(0)=-\frac{l}{k}$, and we let $\p\in\Der (A)$ be a
homogeneous locally nilpotent derivation with $e=\deg \p>0$. Then
there exists a unitary polynomial $Q\in\C[t]$ with $Q(0)\neq 0$
and $\div (t^lQ(t))=-kD_-$ such that, if $A'= A_{k,P}$ is the
normalization of
\be
 \quad\label{Bkp} B_{k,P}=\C[u,v, s]/ \left(
u^kv-P(s)\right)\,, \quad\mbox{where}\quad
P(s)=Q(s^d)s^{ke'+dl}\,,
\ee
then the group $\Z_d=\langle\zeta\rangle$ acts on $B_{k,P}$
and also on $A'$ via
\be
\label{opac} \zeta . (u,v,s)=(\zeta^{e'}u,\,v,\,\zeta s)\,,
\ee
so that $A\cong A^{\prime \Z_d}$. Furthermore, $ee'\equiv 1
\mod d$ and $\p=cu^e\frac{\p}{\p s}\vert A$ for some constant $c\in\C^*$.
\elem

With this result we can complete the proof of Theorem \ref{MT} as
follows. We may assume that $A=A_0[D_+,D_-]$ with $A_0=\C[t]$ and
$(D_+,D_-)$ as in Lemma \ref{L1}. With $k:=\lcm(d,m)$ let us
write $k=mm'=dd'$ and $l=-e'd',$ so that
$$D_+=-\frac{e'}{d}
[0]=\frac{l}{k} [0],\qquad D_-=\frac{e'}{d} [0] -\frac{1}{m}
[1]=-\frac{l}{k} [0]-\frac{m'}{k} [1]\,.$$ Thus Lemma \ref{L2} can
be applied in our setting with $Q=(t-1)^{m'}$. By this lemma,
$A=A'^{\Z_d}$, where $A'$ is the normalization of
$$
B=\C[u,v,s]/(u^kv-(s^d-1)^{m'}),
$$
with the action of $\Z_d$ as in (\ref{opac}) and with the
$\C^*$-action $\lambda.(u,v,s)=(\lambda u, \lambda^{-k} v,s)$.

The element $w=\frac{s^d-1}{u^m}\in\Frac (B)$ satisfies $w^{m'}=v$
and so is integral over $B$, hence
$$
A'\cong\C[u,w,s]/(u^mw-(s^d-1)).
$$
Because of (\ref{opac}) we have $\zeta.w=\zeta^{-me'}w$. Thus
after applying an automorphism $\zeta\mapsto \zeta^{e'}$ of
$\Z_d$, both the $\Z_d$-action and the $\C^*$-action on $\tX=\Spec
A'\subseteq \A^3_\C=\Spec \C[u,w,s]\cong \Spec \C[x,y,z]$ have
the claimed form $$\zeta.(u,w,s)=(\zeta u,\zeta^{-m} w, \zeta^e
s)\qquad\mbox{respectively,}\qquad \lambda.(u,w,s)=(\lambda
u,\lambda^{-m} w, s)\,.$$ This proves the theorem. \qed

\end{document}